\documentclass{amsart}
\usepackage[all]{xy}
\usepackage{amsmath,amsfonts,epsfig,amscd}

\begin{document}
\bibliographystyle{plain}

\newtheorem{thm}{Theorem}
\newtheorem{lem}[thm]{Lemma}
\newtheorem{prop}[thm]{Proposition}
\newtheorem{cor}[thm]{Corollary}
\newtheorem{conj}[thm]{Conjecture}
\newtheorem{mainlem}[thm]{Main Lemma}
\newtheorem{defn}[thm]{Definition}
\newtheorem{rmk}[thm]{Remark}

\newtheorem*{namedtheorem}{\theoremname}
\newcommand{\theoremname}{testing}
\newenvironment{named}[1]{\renewcommand{\theoremname}{#1}\begin{namedtheorem}}{\end{namedtheorem}}


\def\square{\hfill${\vcenter{\vbox{\hrule height.4pt \hbox{\vrule width.4pt
height7pt \kern7pt \vrule width.4pt} \hrule height.4pt}}}$}

\newenvironment{pf}{{\it Proof:}\quad}{\square \vskip 12pt}
\newcommand{\dt}{\ensuremath{\text{det}}}
\newcommand{\U}{\ensuremath{\widetilde}}
\newcommand{\Hn}{\ensuremath{\mathbb{H}^3}}
\newcommand{\h}{{\text{hyp}}}
\newcommand{\acyl}{{\text{acyl}}}
\newcommand{\inc}{{\text{inc}}}
\newcommand{\tp}{{\text{top}}}
\newcommand{\V}{\ensuremath{\text{Vol}}}
\newcommand{\Hess}{\ensuremath{{\text{Hess} \ }}}
\newcommand{\PM}{\ensuremath{\mathcal{PML}(S)}}
\newcommand{\ilim}{\ensuremath{\underleftarrow{\lim}}}
\newcommand{\Cech}{\ensuremath{\check{\text{C}}\text{ech\ }}}

\title[Higher signatures of Mod$(S)$]{The Novikov conjecture for mapping class groups as a corollary of {H}amenst{\"a}dt's theorem}
\author{Peter A. Storm}

\date{December 16th, 2004}



\thanks{The author receives partial funding from an NSF Postdoctoral Research Fellowship.}


\maketitle

Let Mod$(S)$ denote the mapping class group of a finite type surface.  This short note combines theorems of Kato and Hamenst{\"a}dt to show

\begin{cor} \label{main cor}
The higher signatures of Mod$(S)$ are oriented homotopy invariants.  In other words, the Novikov conjecture is true for Mod$(S)$.
\end{cor}

The corollary virtually follows immediately from the following theorems.  The first is due to Kato.

\begin{thm} \label{Kato thm} \cite[Thm.0.1]{Kato}
Let $\Gamma$ be a torsion free finitely generated group which admits a proper combing of bounded multiplicity.  Then the higher signatures of $\Gamma$ are oriented homotopy invariants.
\end{thm}

{\noindent}We will not define a ``proper combing of bounded multiplicity".  Instead we simply note that any ``quasi-geodesically bicombable group can admit a structure of proper combing [\textit{sic}] of strictly bounded multiplicity \cite[Ex.2.1,pg.68]{Kato}."  The second theorem is due to Hamenst{\"a}dt.

\begin{thm} \label{Ham thm} \cite[Sec.6]{Ham}
The mapping class group Mod$(S)$ admits a quasigeodesic bicombing.
\end{thm}

{\noindent}In Section 6 of \cite{Ham}, Hamenst{\"a}dt builds a quasigeodesic bicombing of the train track complex (notated $\mathcal{T}T$), which she shows to be quasi-isometric to Mod$(S)$.  

\begin{rmk}
The proof of Corollary \ref{main cor} rests entirely on Kato's Thorem \ref{Kato thm} and Hamenst{\"a}dt's Theorem \ref{Ham thm}.  It is hoped that this note will bring further attention to these two important results.
\end{rmk}

We now recall the statement of the Novikov conjecture.  Let $M$ be a closed oriented smooth $n$-manifold with fundamental class $[ M] \in H_n (M ; \mathbb{Q})$.  Let $L_M \in H(M ; \mathbb{Q})$ denote the Hirzebruch $L$-class of $M$, which is defined as a (fixed) power series in the Pontryagin classes of $M$.  Let $\Gamma$ be a discrete group with an Eilenberg-Maclane space $B\Gamma$.  For a $u \in H(B\Gamma ; \mathbb{Q})$ and a continuous $f:M \longrightarrow B\Gamma$ define the higher signature
$$ \langle L_M \cup f^* u , [M] \rangle \in \mathbb{Q}.$$
The higher signature defined by $u$ and $f$ is an oriented homotopy invariant if: for any closed oriented smooth manifold $N$ with fundamental class $[N]$ and a homotopy equivalence $h: N \longrightarrow M$ taking $[N]$ to $[M]$ we have the equality
$$\langle L_M \cup f^* u , [M] \rangle = \langle L_N \cup (f \circ h)^* u , [N] \rangle. \qquad (\dagger)$$
For more information on the Novikov conjecture see \cite{Novikov}.

A finite presentation for Mod$(S)$ was given by Hatcher-Thurston \cite{HT}.  Note also the following well known fact.

\begin{prop} \cite[Ch.5.5]{K} \label{torsion free} 
The mapping class group contains a finite index normal torsion free subgroup $\Lambda$.
\end{prop}

{\noindent}Combining Theorem \ref{Kato thm}, Theorem \ref{Ham thm}, and the above proposition establishes the Novikov conjecture for $\Lambda$.  It remains only to extend the result to all of Mod$(S)$.

Define $\Gamma := \text{Mod}(S)$ and let $M, N, f, h$, and $u$ be as above.  We must establish equation $(\dagger)$.  Let $\U{M}$, $\U{N}$, $\U{f}$, and $\U{h}$ denote the lifts of $M$,$N$, $f$, and $h$ respectively to the finite covering spaces corresponding to $(f_*)^{-1}(\Lambda) \le \pi_1 (M)$, where $f_* : \pi_1 (M) \longrightarrow \Gamma$ is the induced map on $\pi_1$.  We will use $\pi$ to denote either of the covering maps $\U{M} \longrightarrow M$, $\U{N} \longrightarrow N$, or $B \Lambda \longrightarrow B\Gamma$.  Applying the Novikov conjecture for $\Lambda$ to these lifts yields the equation
	$$\langle L_{\U{M}} \cup \U{f}^* (\pi^* u) , [\U{M}] \rangle = 
	    \langle L_{\U{N}} \cup (\U{f} \circ \U{h})^* (\pi^* u) , [\U{N}] \rangle.$$
The naturality of the pairing implies
\begin{eqnarray*}
\langle \pi^* L_M \cup \U{f}^* (\pi^* u) , [\U{M}] \rangle & = &
 \langle \pi^* L_M \cup \pi^* (f^* u) , [\U{M}] \rangle \\
 & = & \langle L_M \cup f^* u , \pi_* [\U{M}] \rangle = k \cdot \langle L_M \cup f^* u , [M] \rangle,
\end{eqnarray*}
where $k$ is the index of $(f_*)^{-1} (\Lambda)$ in $\pi_1 (M)$.  Similarly
\begin{eqnarray*}
\langle \pi^* L_N \cup (\U{f} \circ \U{h})^* (\pi^* u) , [\U{N}] \rangle & = &
 \langle \pi^* L_N \cup \pi^* ((f \circ h)^* u)  , [\U{N}] \rangle \\
 & = & \langle L_N \cup (f \circ h)^* u  , \pi_* [\U{N}] \rangle \\
 & = & k \cdot \langle L_N \cup (f \circ h)^* u , [N] \rangle.
\end{eqnarray*}
By naturality $L_{\U{M}} = \pi^* L_M$ and $L_{\U{N}} = \pi^* L_N$.  Therefore
\begin{eqnarray*}
k \cdot \langle L_M \cup f^* u , [M] \rangle & = & \langle L_{\U{M}} \cup \U{f}^* (\pi^* u), [\U{M}] \rangle \\
	   & = & \langle L_{\U{N}} \cup (\U{f} \circ \U{h})^* (\pi^* u), [\U{N}] \\
	   & = & k \cdot \langle L_N \cup (f \circ h)^* u , [N] \rangle.
\end{eqnarray*}
This establishes equation $(\dagger)$ and completes the proof.

\vskip 6pt
{\noindent}\textbf{Acknowledgement:}  This note is the product of a conversation between the author and Chris Connell.

\bibliography{storm.biblio} 
\end{document}